\documentclass[english]{article}
\usepackage[T1]{fontenc}
\usepackage[latin9]{inputenc}
\usepackage{amsmath}
\usepackage{amssymb}
\usepackage{graphicx}

\makeatletter
\renewcommand{\Re}{\operatorname{Re}}

\usepackage{cite}
\date{Created November 23, 2006; updated \today}

\makeatother

\usepackage{babel}
\begin{document}

\title{Saddle-point integration\\
of $C_{\infty}$ ``bump'' functions}

\author{Steven G. Johnson, MIT Applied Mathematics}
\maketitle
\begin{abstract}
This technical note describes the application of saddle-point integration
to the asymptotic Fourier analysis of the well-known $C_{\infty}$
``bump'' function $e^{-(1-x^{2})^{-1}}\mathbf{1}_{(-1,1)}$, deriving
both the asymptotic decay rate $k^{-3/4}e^{-\sqrt{k}}$ of the Fourier
transform $F(k)$ and the exact coefficient. The result is checked
against brute-force numerical integration and is extended to generalizations
of this bump function.
\end{abstract}

\section{Background}

In our work involving optimization of slow ``taper'' transitions
for waveguide couplers and also for our work on exponential absorbing
boundary layers in wave equations~\cite{OskooiZh08,OskooiMu12},
we encountered a problem that could be expressed in terms of the Fourier
transform of smooth functions with compact support. In particular,
we needed to analyze the asymptotic rate of decay of the Fourier coefficients.
This was straightforward for functions with simple discontinuities
(via integration by parts~\cite{Katznelson68,MeadDe73,Boyd01}) or
with simple poles (via contour integration~\cite{Elliott64,Boyd01}).
However, it turned out that some of the most interesting functions
were $C_{\infty}$ (infinitely differentiable) functions with compact
support, and these functions have essential singularities that resist
those methods. In this note, we describe how the asymptotic Fourier
transforms of such functions can be analyzed with the help of saddle-point
integration \cite{Cheng06}.\footnote{We are indebted to our colleagues, Prof. Martin Bazant and especially
Prof. Hung Cheng at MIT, for their helpful suggestions in this matter.} A similar approach has been applied to various functions in the related
context of Chebyshev polynomial series~\cite{ElliottSz65,Miller66,Boyd82}.

In particular, we will look at $C_{\infty}$ functions on $\mathbb{R}$
with compact support $[-1,1]$. The canonical example of such a function
is the symmetric ``bump:''
\[
f(x)=\left\{ \begin{array}{c}
e^{-\frac{1}{1-x^{2}}}\quad x\in(-1,1)\\
0\quad\textrm{otherwise}
\end{array}\right..
\]
 The actual functions we are interested may be more complicated than
this {[}e.g., $f(x)$ multiplied by some analytic function{]}, but
their analysis is similar to that of $f(x)$: the key point is that
we have essential singularities at $x=-1$ and $x=+1$, and these
singularities determine the asymptotic behavior of the Fourier transform
and similar integrals.

We should also note that the space of $C_{\infty}$ functions with
compact support plays an important role in the theory of generalized
functions (distributions), where they are the ``test functions''
that are the domain of the distributions. In this context, it has
been proven that the Fourier transform of any such test function is
an entire function (analytic over the whole complex plane) and diverges
at most exponentially fast off the real axis \cite{Gelfand64}. From
the fact that the functions are infinitely differentiable, it also
immediately follows that their Fourier transforms go to zero along
the real axis faster than the inverse of any polynomial \cite{Katznelson68,MeadDe73}.
Here, however, we want to know precisely how fast the Fourier transform
decays, and with what coefficient.

\section{Saddle-point Fourier integration of $f(x)$}

We wish to compute the asymptotic behavior of the Fourier transform
\[
F(k)=\int_{-\infty}^{\infty}f(x)e^{ikx}dx=2\Re\int_{0}^{1}e^{ikx-\frac{1}{1-x^{2}}}dx
\]
for $|\Re k|\gg1$. (Without loss of generality, we can restrict ourselves
to real $k\geq0$.) To do this, we will employ a saddle-point integration:
we will look for the $x$ at which the exponent is stationary, and
expand the exponent approximately around this point. For large $k$,
this saddle point of the stationary exponent will dominate the integral,
and this will give us the asymptotic behavior. It will turn out that
the saddle point occurs for a complex $x$, however, so we will need
to deform the integration contour within the complex plane to exploit
this approach.

It is clear (and this will also be justified \emph{a posteriori})
that the biggest contributions must come near the singular point $x=1$.
Exactly \emph{at} the singular point, however, the integrand is zero,
so (perhaps) contrary to our intuition the endpoints \emph{per se}
are not important. Because we are worrying about points near the endpoint,
however, it is convenient to change variables $t=1-x$ and write 
\[
F(k)=2\Re\int_{0}^{1}e^{g(t)}dt
\]
with 
\[
g(t)=ik-ikt-\frac{1}{(2-t)t}\approx ik-ikt-\frac{1}{2t}-\frac{1}{4}+O(t),
\]
 where the last approximation is for $|t|\ll1$, which is valid (it
turns out) at our saddle point. The saddle point is the $t=t_{0}$
where $g'(t)=0=-ik+1/2t^{2}$, which by inspection is\footnote{$g'$ also vanishes at $t=-\sqrt{1/2ik}$, but we cannot deform our
integration contour to a point with negative $\Re t$: for $\Re t<0$,
our integrand $\sim e^{-1/2t}$ blows up.} $t_{0}=\sqrt{1/2ik}$. Note that for large $k\gg1$ we obtain $|t_{0}|\ll1$,
justifying our approimation above. Now, we write 
\begin{eqnarray*}
g(t) & \approx & g(t_{0})+\frac{g''(t_{0})}{2}(t-t_{0})^{2}\\
 & = & ik-i\sqrt{\frac{-ik}{2}}-\sqrt{\frac{ik}{2}}-\frac{1}{4}+\frac{-(\sqrt{2ik})^{3}}{2}\left(t-\sqrt{\frac{1}{2ik}}\right)^{2}\\
 & = & ik-\frac{1}{4}-\sqrt{2ik}+\sqrt{-2i}k^{3/2}\left(t-\sqrt{\frac{1}{2ik}}\right)^{2}.
\end{eqnarray*}
 To actually do this integral, we need to deform our integration contour
in the complex plane to lie along a line $t=u/\sqrt{i}$ for real
$u$ near $u=0$, so that the saddle point $t=t_{0}$ lies along our
integration path (at $u=\sqrt{1/2k}$).\footnote{Although $t=u/\sqrt{i}$ is not the path of \emph{steepest} descent
around $t_{0}$, as would be prescribed by the usual saddle-point
method (a.k.a. the method of ``steepest descent''), it is at least
a path of descent (the integrand is decaying along that path) \cite{Cheng06}.
This is sufficient for us to apply our Gaussian integral approximation.} (This deformation is not a problem since we don't have any singularities
except at $t=0$ and $t=2$, and the integrand vanishes as $t\rightarrow0$
for $\Re t>0$.) In this case, our integral becomes (approximately)
a Gaussian integral, since: 
\[
g(u/\sqrt{i})=ik-\frac{1}{4}-\sqrt{2ik}-\sqrt{2i}k^{3/2}\left(u-\sqrt{\frac{1}{2k}}\right)^{2}.
\]
 Recall that the integral of $\int_{-\infty}^{\infty}e^{-au^{2}}du=\sqrt{\pi/a}$
as long as $\Re a>0$, which is true here. Note also that, in the
limit as $k$ becomes large, the integrand becomes zero except close
to $u=\sqrt{1/2k}$, so we can neglect the rest of the contour and
treat the integral over $u$ as going from $-\infty$ to $\infty$.
(Thankfully, the width of the Gaussian $\Delta u\sim k^{-3/4}$ goes
to zero faster than the location of the maximum $u_{0}\sim k^{-1/2}$,
so we don't have to worry about the $u=0$ origin.) Also note that
the change of variables from $t$ to $u$ gives us a $\sqrt{-i}$
Jacobian factor. Thus, when all is said and done, we obtain the exact
asymptotic form of the Fourier integral for $k\gg1$:
\[
F(k)\approx2\Re\left[\sqrt{\frac{-i\pi}{\sqrt{2i}k^{3/2}}}e^{ik-\frac{1}{4}-\sqrt{2ik}}\right].
\]
Since $\sqrt{2i}=1+i$, this means that the Fourier coefficients decay
proportional to $k^{-3/4}e^{-\sqrt{k}}$, which is consistent with
the expected faster-than-polynomial decay.

\begin{figure}
\begin{centering}
\includegraphics[width=1\columnwidth]{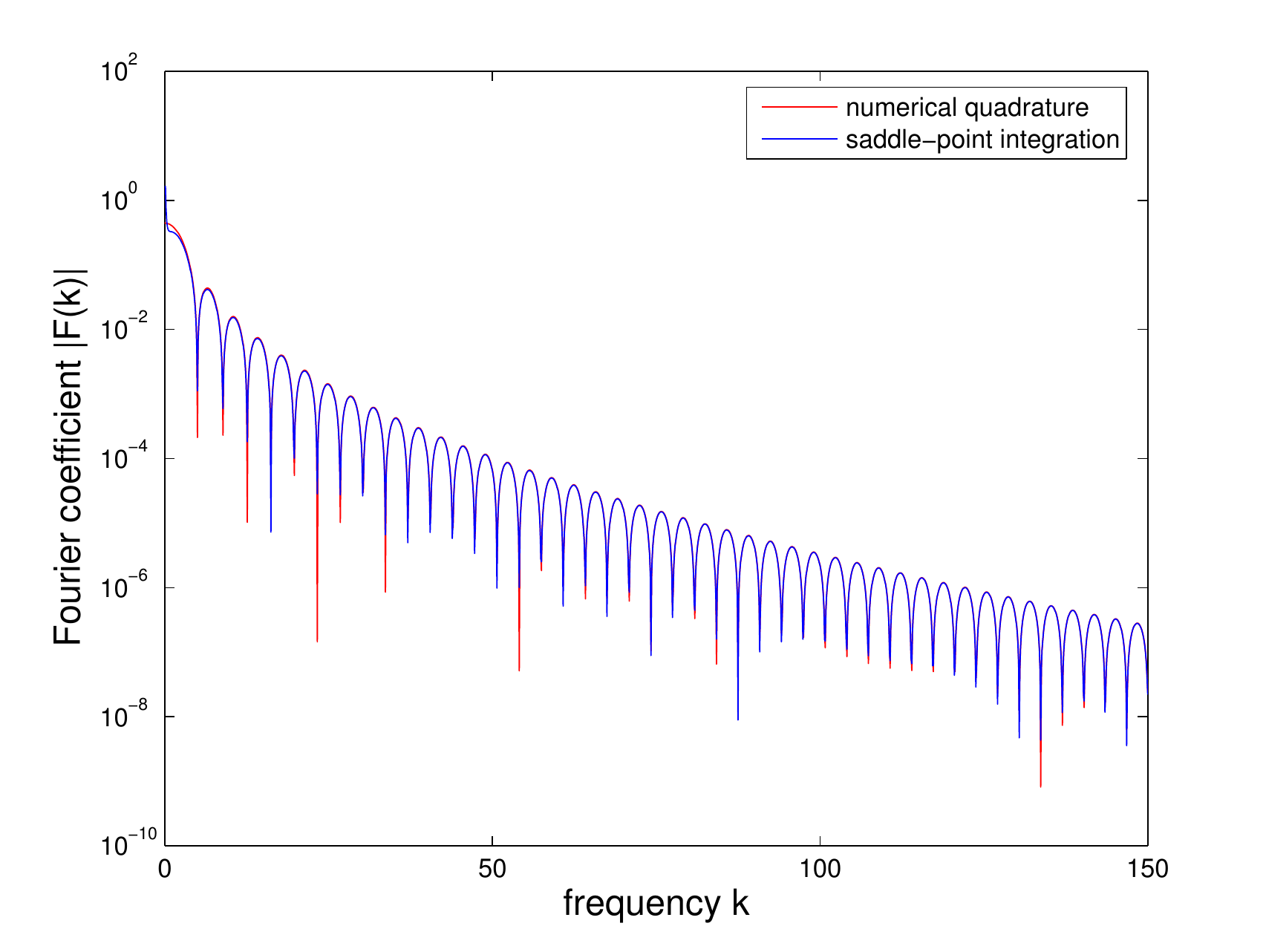}
\par\end{centering}

\protect\caption{\label{fig:saddle-bump-quadrature}``Exact'' numerical quadrature
of $|F(k)|$ (red line) and asymptotic saddle-point integration (blue
line) for the Fourier transform of the $C_{\infty}$ bump function
$f(x)$. (They match so well that the curves are difficult to distinguish.)}
\end{figure}
To check this result, we can compare the above formula with an exact
numerical evaluation of the Fourier integral. Numerical integration
was performed using adaptive Clenshaw-Curtis quadrature, specialized
for a $\cos(kx)$ oscillatory weight function, from the GNU Scientific
Library (adapted from QUADPACK). The resulting $|F(k)|$ for $0\leq k\leq150$
is plotted in figure~\ref{fig:saddle-bump-quadrature}. As can be
seen from the figure, the asymptotic approximation is an excellent
match for the exact result, with errors under 10\% for frequencies
$k$ as small as $4$.

\section{Generalized bumps}

As a further warm-up for the more complicated problems we may want
to solve later, let's look at the asymptotic Fourier transform of
a generalization of our bump function $f(x)$:
\[
f_{\alpha,\beta}(x)=\left\{ \begin{array}{c}
e^{-\frac{\beta}{(1-x^{2})^{\alpha-1}}}\quad x\in(-1,1)\\
0\quad\textrm{otherwise}
\end{array}\right.,
\]
so that $f(x)=f_{2,1}$. Again, we look at the Fourier transform with
changed variables $x=1-t$, i.e. $\int_{0}^{1}f_{\alpha,\beta}(1-t)e^{ik-ikt}dt$,
and expand the exponent $g(t)$ around $t=0$: 
\[
g(t)=ik-ikt-\frac{\beta}{(2-t)^{\alpha-1}t^{\alpha-1}}\approx ik-ikt-\frac{\beta}{(2t)^{\alpha-1}}+O(t^{2-\alpha}).
\]
 Our derivative is then $g'(t)=0\approx-ik+\beta(\alpha-1)/2^{\alpha-1}t^{\alpha}$
at 
\[
t_{0}=[2^{1-\alpha}\beta(\alpha-1)/ik]^{1/\alpha}.
\]
 That is, $t_{0}\sim k^{-1/\alpha}$, and therefore $g(t_{0})\sim k^{(\alpha-1)/\alpha}$
asymptotically.

The second derivative is $g''(t_{0})=-2^{1-\alpha}\beta\alpha(\alpha-1)/t_{0}^{\alpha+1}=-i^{(\alpha+1)/\alpha}2Ak^{(\alpha+1)/\alpha}$
where 
\[
A=\alpha[2\beta(\alpha-1)]^{-1/\alpha}.
\]
 Again, we'll choose a contour $t=u/i^{1/\alpha}$, in which case
we find: 
\[
g(t)\approx g(t_{0})-i^{(\alpha-1)/\alpha}Ak^{(\alpha+1)/\alpha}(u-u_{0})^{2},
\]
which is a path of descent so we can perform a Gaussian integral.
The final answer for the integral, including the Jacobian factor for
$dt=du/i^{1/\alpha}$, is then 
\[
F(k)\approx2\Re\left[\sqrt{\frac{\pi}{(ik)^{(\alpha+1)/\alpha}A}}e^{ik-ikt_{0}-\frac{\beta}{[(2-t_{0})t_{0}]^{\alpha-1}}}\right],
\]
with $t_{0}$ and $A$ given above. Notice that for our exponent $g(t_{0})$
we use the exact form of $g(t)$ and not its approximation for small
$t$. The reason is that, for $\alpha>2$, terms of $O(t^{2-\alpha})$
do not go to zero, and therefore make a non-negligible multiplicative
contribution to the amplitude $F(k)$ even though they do not affect
the saddle-point integration.

Numerical tests seem to confirm the accuracy of this $F(k)$ formula,
although for $\alpha>3$ we start to have more difficulty with the
numerical quadrature. As one might have expected, increasing either
$\alpha$ or $\beta$ makes $F(k)$ decay more rapidly. To perform
such comparisons more carefully, we should typically normalize by
$\int_{0}^{1}f_{\alpha,\beta}(x)dx$ or similar, to ignore effects
due simply to the fact that the integrand is getting smaller overall.

\bibliographystyle{ieeetr}
\bibliography{photon}

\end{document}